\documentclass[11pt,reqno]{amsart}
\usepackage{amsmath,amssymb,geometry,color,tikz-cd}
\usepackage[backref,pagebackref,pdftex,hyperindex]{hyperref}
\usepackage[alphabetic]{amsrefs}
\usepackage{amssymb}
\usepackage{bbm}
\usepackage{enumitem}
\usepackage{upgreek}

\newtheorem{proposition}{Proposition}[section]
\newtheorem{theorem}[proposition]{Theorem}
\newtheorem{lemma}[proposition]{Lemma}

\theoremstyle{definition}
\newtheorem{remark}[proposition]{Remark}
\newtheorem{definition}[proposition]{Definition}
\newtheorem{example}[proposition]{Example}

\title{Log K-stability of GIT-stable divisors on Fano varieties}
\author{Chuyu Zhou}
\address{\'Ecole Polytechnique F\'ed\'erale de Lausanne (EPFL), MA C3 615, Station 8, 1015 Lausanne, Switzerland}
\email{chuyuzhou1@gmail.com}
\date{} 


\newcommand{\Fut}{{\rm{Fut}}}

\newcommand{\ord}{{\rm {ord}}}
\newcommand{\tc}{{\rm {tc}}}
\newcommand{\vol}{{\rm {vol}}}

\newcommand{\dt}{{\rm {dt}}}

\newcommand{\Hilb}{{\rm {Hilb}}}
\newcommand{\Chow}{{\rm {Chow}}}
\newcommand{\PGL}{{\rm {PGL}}}

\newcommand{\Aut}{{\rm {Aut}}}

\newcommand{\bA}{\mathbb{A}}

\newcommand{\bC}{\mathbb{C}}

\newcommand{\bN}{\mathbb{N}}

\newcommand{\bP}{\mathbb{P}}
\newcommand{\bQ}{\mathbb{Q}}

\newcommand{\mD}{\mathcal{D}}

\newcommand{\mF}{\mathcal{F}}

\newcommand{\mL}{\mathcal{L}}
\newcommand{\mM}{\mathcal{M}}

\newcommand{\mO}{\mathcal{O}}

\newcommand{\mX}{\mathcal{X}}
\newcommand{\mY}{\mathcal{Y}}

\newcommand{\tD}{\tilde{D}}

\begin{document}

\begin{abstract}
For a given K-polystable Fano variety $X$ and a  natural number $l$ such that $(X, \frac{1}{l}B)$ is log canonical for some $B\in |-lK_X|$, we show that there exists a rational number $0<c_1<1$ depending only on $X$ and $l$, such that $D\in |-lK_X|$ is GIT-(semi/poly)stable under the action of $\Aut(X)$ if and only if the log pair $(X, \frac{\epsilon}{l}D)$ is K-(semi/poly)stable for any rational $0<\epsilon<c_1$.
\end{abstract}

\maketitle
\tableofcontents

\section{Introduction}\label{sec:1}

In the past few years, people have made tremendous progress on the construction of K-moduli space of Fano varieties, e.g \cite{Jiang20, BLX22, Xu20, BX19, ABHLX20, XZ20b, LXZ21}, and it turns out that there exists a projective scheme (as a good moduli space) parametrizing Fano varieties with K-stability and fixed invariants.  In particular, the K-moduli space is proper. This is first known for smoothable Fano varieties with the help of analytic tools, e.g \cite{CDS15a,CDS15b,CDS15c, Tian15, TW19, LWX19, ADL19} ect., and later proven by \cite{LXZ21} for Fano varieties without smoothable condition.
With properness in hand, moduli continuity method has been widely applied in the literature to construct explicit K-moduli spaces of some special kinds of Fano varieties, e.g.  K-moduli for Del Pezzo surfaces (\cite{OSS16}), for cubic 3-folds (\cite{LX19}) and cubic 4-folds (\cite{Liu20}). It is well known that GIT-stability and K-stability, although not the same, have close relationship via CM-line bundle. More precisely, the generalized Futaki invariant of an one parameter subgroup can be identified with the corresponding GIT-weight via the CM-line bundle on the base (up to a positive multiple), thus in many cases one can construct a morphism from a K-moduli space to a GIT-moduli space. The moduli continuity method aims to establish an isomorphism between these two spaces, which becomes a powerful way to confirm K-stability of explicit Fano varieties. 

In order to state the main theorem, we first fix some notation. Let $X$ be a K-polystable Fano variety and $l$ a positive natural number such that $(X, \frac{1}{l}B)$ is log canonical for some $B\in |-lK_X|$. As $\Aut(X)$ is reductive (e.g \cite{ABHLX20, CDS15a, CDS15b, CDS15c, Tian15, Matsushima57}), it is natural to determine GIT-stability of elements in $|-lK_X|$ under the action of $\Aut(X)$. For a rational number $0<\epsilon<1$, we denote by $\mM_{X,l,\epsilon}^K$ the Artin stack which parametrizes all K-semistable log Fano pairs of the form $(Y,\frac{\epsilon}{l}\tD)$, where $(Y,\tD)$ is the degeneration of $(X,|-lK_X|)$ (see Section \ref{sec:4}). We also write $\mM^{GIT}_{X,l}:=[(|-lK_X|)^{ss}/\Aut(X)]$. Denote by $M^K_{X,l,\epsilon}$ (resp. $M^{GIT}_{X,l}$) the good moduli space of $\mM^K_{X,l,\epsilon}$ (resp. $\mM^{GIT}_{X,l}$). It is recognized that the GIT-stability of $D\in |-lK_X|$ is related to the K-stability of the log Fano pair $(X, \frac{\epsilon}{l}D)$ for small $\epsilon$. For example, in the case $X=\bP^n$, the work \cite{GMGS21} establishes such a correspondence between GIT-stability and K-stability for $n=1,2,3$, and the work \cite[Theorem 1.4]{ADL19} establishes it for every $n$.
In this paper, we apply moduli continuity method to prove the following theorem, which is conjectured in \cite[Conjecture 1.3]{GMGS21}. 

\begin{theorem}\label{main thm}
Let X be a K-polystable Fano variety and $l$ a positive integer such that $(X, \frac{1}{l}B)$ is log canonical for some $B\in |-lK_X|$.  Then there exists a rational number $0<c_1<1$ depending only on  $X$ and $l$ such that the following two statements are equivalent:
\begin{enumerate}
\item $D\in |-lK_X|$ is GIT-(semi/poly)stable under the action of $\Aut(X)$,
\item the pair $(X, \frac{\epsilon}{l}D)$ is K-(semi/poly)stable for any rational $0<\epsilon< c_1$.
\end{enumerate}
As a consequence, we have an isomorphism $\phi_\epsilon: \mM^K_{X,l,\epsilon} \to \mM^{GIT}_{X,l}$ which induces an isomorphism on the moduli spaces $\phi_\epsilon': M^K_{X,l,\epsilon} \to M^{GIT}_{X,l}$ for any rational $0<\epsilon< c_1$.
\end{theorem}

The above result generalizes \cite[Theorem 1.4]{ADL19} to any K-polystable Fano variety. 
The key ingredient for the proof of the theorem is to confirm that, the K-semistable degeneration of $(X,\frac{\epsilon}{l}D)$ for small $\epsilon$ and $D\in |-lK_X|$ preserves the ambient space. The approach in \cite{ADL19} for this ingredient is based on the fact that any KE Fano variety with maximal volume is isomorphic to projective space (\cite{Liu18, Fuj18}), while our approach is to reduce the problem to the fact that,
if a K-polystable Fano variety admits a K-semistable degeneration,  then the degeneration is still K-polystable (see Theorem \ref{key}).
In fact, if the degeneration is obtained by a test configuration, then the fact is well known by \cite[Section 3]{LWX21} (see also Definition \ref{def ks}).   

We also pose the condition that $(X, \frac{1}{l}B)$ is log canonical for some $B\in |-lK_X|$, and it makes sure that the stack $\mM^K_{X, l, \epsilon}$ is non-empty (see Lemma \ref{non-empty}).  The non-emptiness indeed plays a role in the final proof of Theorem \ref{main thm} (see the second paragraph of the proof of Theorem \ref{thm: main thm}). If we remove this condition, it is interesting to find an example where $\mM^K_{X, l, \epsilon}$ is empty while $\mM^{GIT}_{X, l}$ not.

\noindent
\subsection*{Acknowledgement}
The author would like to thank Chen Jiang, Yuchen Liu, Ziquan Zhuang for helpful discussions and beneficial comments. The author is supported by grant European Research Council (ERC-804334).

\section{Preliminaries}\label{sec:2}

In this section, we provide necessary preliminaries. We work over complex number field $\bC$. We say that $(X,\Delta)$ is a log pair if $X$ is a projective normal variety and $\Delta$ is an effective $\bQ$-divisor on $X$ such that $K_X+\Delta$ is $\bQ$-Cartier. We say a log pair $(X,\Delta)$ is log Fano if it admits klt singularities and $-K_X-\Delta$ is ample. If $\Delta=0$, we just say a log Fano pair $(X,\Delta)$ is a Fano variety. For the concepts of singularities in birational geometry such as klt singularities, we refer to \cite{KM98, Kollar13}.

\subsection{K-stability of Fano varieties}\label{subsec:ks}

Let $(X,\Delta)$ be a log Fano pair of dimension $n$, we denote $L:=-K_X-\Delta$, which is an ample $\bQ$-line bundle.

\begin{definition}
We say that a triple $\pi: (\mX,\Delta_{\tc};\mL)\to \bA^1$ is a test configuration of $(X,\Delta;L)$ if the following conditions are satisfied:
\begin{enumerate}
\item $\pi$ is a flat projective morphism from a normal variety $\mX$ and $\Delta_{\tc}\subset \mX$ is a $\bQ$-divisor flat over $\bA^1$,
\item $\mL$ is a relatively ample $\bQ$-line bundle on $\mX$ with a $\bC^*$-action induced by the natural multiplication on $\bA^1$,
\item $(\mX^*,\Delta_{\tc}^*;\mL^*)$ is $\bC^*$-equivariantly isomorphic to $(X\times \bC^*, \Delta\times \bC^*; L\times \bC^*)$, where $\mX^*:=\mX\setminus \mX_0$.
\end{enumerate}
The compactification of the test configuration is denoted by $(\bar{\mX}, \bar{\Delta}_{\tc};\bar{\mL})\to \bP^1$, which is obtained by gluing $(\mX,\Delta_{\tc})$ and $(X\times (\bP^1\setminus 0), \Delta\times (\bP^1\setminus 0))$ along $(X\times \bC^*, \Delta\times \bC^*)$.
\end{definition}

\begin{definition}\label{def ks}
Let $(\mX,\Delta_{\tc}; \mL)\to \bA^1$ be a test configuration of $(X,\Delta;L)$, then the generalized Futaki invariant of this test configuration is defined as follows:
$$\Fut(\mX,\Delta_{\tc};\mL):=\frac{n\bar{\mL}^{n+1}}{(n+1)(-K_X-\Delta)^n} +\frac{\bar{\mL}^n(K_{\bar{\mX}/\bP^1}+\bar{\Delta}_{\tc})}{(-K_X-\Delta)^n}.$$
We say that $(X,\Delta)$ is K-semistable if $\Fut(\mX,\Delta_{\tc};\mL)\geq 0$ for every test configuration. We say that $(X,\Delta)$ is K-polystable if it is K-semistable and for any test configuration of $(X,\Delta;L)$ whose central fiber is K-semistable, we have an isomorphism between $(X,\Delta)$ and the central fiber\footnote{This is not the original definition of K-polystability, but equivalent to the original one by \cite{LWX21}.}.
\end{definition}

\subsection{Valuative criterion}\label{sec value}

Let $(X,\Delta)$ be a log Fano pair of dimension $n$. We say $E$ is a prime divisor over $X$ if there is a proper birational morphism from a normal variety $\phi: Y\to X$ such that $E$ is a prime divisor on $Y$. We define
$$A_{X,\Delta}(E):=\ord_E(K_Y-\phi^*(K_X+\Delta))+1 ,$$
$$S_{X,\Delta}(E):=\frac{1}{(-K_X-\Delta)^n}\int_0^{\infty}\vol(-\phi^*(K_X+\Delta)-tE)\dt.$$

\begin{definition}
The beta-invariant of a prime divisor $E$ over a log Fano pair $(X,\Delta)$ is defined as follows:
$$\beta_{X,\Delta}(E):=A_{X,\Delta}(E)-S_{X,\Delta}(E).$$
The delta invariant of a log Fano pair $(X,\Delta)$ is defined as follows:
$$\delta(X,\Delta):=\inf_E\frac{A_{X,\Delta}(E)}{S_{X,\Delta}(E)}, $$
where $E$ runs through all prime divisors over $X$.
\end{definition}

We have the following well-known theorem due to \cite{Fuj19, Li17, FO18, BJ20}.

\begin{theorem}
Let $(X,\Delta)$ be a log Fano pair of dimension $n$, then
\begin{enumerate}
\item $(X,\Delta)$ is K-semistable if and only if $\beta_{X,\Delta}(E)\geq0$ for any prime divisor $E$ over $X$.
\item $(X,\Delta)$ is K-semistable if and only if $\delta(X,\Delta)\geq 1$.\end{enumerate}
\end{theorem}

\subsection{CM-line bundle}

Let $\pi: (\mX,\mD;\mL)\to T$ be a flat family of projective normal varieties of dimension $n$ over a normal base $T$, where $\mD$ is an effective $\bQ$-divisor on $\mX$ whose components are all flat over $T$, and $\mL$ is a relative ample $\bQ$-line bundle. By the work of Mumford-Knudsen(\cite{KM76}) there exist $\bQ$-line bundles $\lambda_i,i=0,1,...,n+1$ and $\tilde{\lambda}_i,i=0,1,...n,$ on $T$ such that we have the following expansions for all sufficiently large $k\in \bN$:
$$\det \pi_*(\mL^k)= \lambda_{n+1}^{\binom{k}{n+1}}\otimes\lambda_n^{\binom{k}{n}}\otimes...\otimes\lambda_1^{\binom{k}{1}}\otimes\lambda_0,$$
$$\det \pi_*(\mL|_\mD^k)= \tilde{\lambda}_{n}^{\binom{k}{n}}\otimes\tilde{\lambda}_{n-1}^{\binom{k}{n-1}}\otimes...\otimes\tilde{\lambda}_0.$$
By Riemann-Roch formula, cf \cite[Appendix]{CP21}, we have
$$c_1(\pi_*\mL^k)=\frac{\pi_*(\mL^{n+1})}{(n+1)!}k^{n+1}+\frac{\pi_*(-K_{\mX/T}\mL^n)}{2n!}k^n+... ,$$
$$c_1({\pi}_*{\mL}|_\mD^k)=\frac{{\pi}_*(\mL^n\mD)}{n!}k^n+... .$$
From above formulas, it's not hard to see
$$\lambda_{n+1}=\pi_*(\mL^{n+1}), \lambda_n=\frac{n}{2}\pi_*(\mL^{n+1})+\frac{1}{2}\pi_*(-K_{\mX/T}\mL^n) \quad and \quad 
\tilde{\lambda}_n=\pi_*(\mL^n\mD). $$
By the flatness of $\pi$ and $\pi_\mD$, we write
$$h^0(\mX_t, k\mL_t)=a_0k^n+a_1k^{n-1}+o(k^{n-1})\quad and \quad  h^0(\mD_t, k{\mL_t}|_{\mD_t})=\tilde{a}_0k^{n-1}+o(k^{n-1}),$$
which do not depend on the choice of $t\in T$. Then we have
$$a_0=\frac{\mL_t^n}{n!}, a_1=\frac{-K_{\mX_t}{\mL_t}^{n-1}}{2(n-1)!} \quad and \quad \tilde{a}_0=\frac{\mL_t^{n-1}\mD_t}{(n-1)!}.$$

\begin{definition}
We define the CM-line bundles for the family $\pi: (\mX,\mD;\mL)\to T$ as follows:
$$\lambda_{CM,(\mX,\mL;\pi)}:=\lambda_{n+1}^{\frac{2a_1}{a_0}+n(n+1)}\otimes\lambda_n^{-2(n+1)}, $$
$$\lambda_{CM,(\mX,\mD,\mL;\pi)}:= \lambda_{n+1}^{\frac{2a_1-\tilde{a}_0}{a_0}+n(n+1)}\otimes\lambda_n^{-2(n+1)}\otimes\tilde{\lambda}_n^{n+1}.$$
\end{definition}

Now we assume $\pi: (\mX,\mD;\mL)\to \bP^1$ to be a compactification test configuration of a log pair $(X,D;L)$ where $L$ is an ample $\bQ$-line bundle on $X$. As $\pi_*\mL^k$ is a $\bC^*$-equivariant vector bundle on $\bP^1$,  we write the total weights to be:
$$w(\det(\pi_*\mL^k))=b_0k^{n+1}+b_1k^n+o(k^n) \quad and \quad w(\det(\pi_*\mL|_\mD^k))=\tilde{b}_0k^n+o(k^n).$$
It is not hard to compute that
\begin{enumerate}
\item $b_0=\frac{\mL^{n+1}}{(n+1)!}=\frac{w(\lambda_{n+1})}{(n+1)!}$,
\item $b_1=\frac{-K_{\mX/\bP^1}\mL^n}{2n!}=\frac{w(\lambda_n)}{n!}-\frac{n(n+1)}{2(n+1)!}w(\lambda_{n+1})$,
\item $\tilde{b}_0=\frac{\mL^n\mD}{n!}=\frac{w(\tilde{\lambda}_n)}{n!}$.
\end{enumerate}

\begin{definition}
Notation as above, we define
\begin{enumerate}
\item The generalized Futaki invariant of $(\mX,\mL)$:
$$\Fut(\mX,\mL):=\frac{2(b_0a_1-b_1a_0)}{a_0^2}=\frac{1}{(n+1)L^n}w(\lambda_{n+1}^{\frac{2a_1}{a_0}+n(n+1)}\otimes\lambda_n^{-2(n+1)}),$$

\item 
The Chow weight of $\pi: (\mX,\mD;\mL)\to T$:
$$\Chow(\mX,\mD;\mL):=\frac{\tilde{b}_0a_0-\tilde{a}_0b_0}{a_0^2}=\frac{1}{(n+1)L^n}w(\lambda_{n+1}^{-\frac{\tilde{a}_0}{a_0}}\otimes\tilde{\lambda}_n^{n+1}),$$

\item The generalized Futaki invariant of $(\mX,\mD;\mL)$:
\begin{align*}
\Fut(\mX,\mD;\mL):=&\Fut(\mX,\mL)+\Chow(\mX,\mD;\mL)\\
=&\frac{1}{(n+1)L^n}w(\lambda_{n+1}^{\frac{2a_1-\tilde{a}_0}{a_0}+n(n+1)}\otimes\lambda_n^{-2(n+1)}\otimes\tilde{\lambda}_n^{n+1}). 
\end{align*}
\end{enumerate}
\end{definition}

\begin{remark}\label{remark on CM}
We have a few remarks for the above definition.
\begin{enumerate}
\item If $\pi: (\mX,\mD;\mL)\to T$ is a compactification test configuration of a log Fano pair, then the generalized Futaki invariants here coincide with the definition in Section \ref{subsec:ks}.
\item In the case of test configurations, we see that the generalized Futaki invariants coincide with the GIT-weights of CM-line bundles up to a multiple. We also note here that CM-line bundles are in fact $\bQ$-line bundles on the base.

\item Suppose $\pi: (\mX,\mD;\mL)\to T$ is a flat family of log Fano varieties such that $\mL\sim_\bQ -K_{\mX/T}-\mD$, then
$$\lambda_{CM,(\mX,\mD,\mL;\pi)}=-\pi_*(\mL^{n+1})\sim_\bQ -\pi_*(-K_{\mX/T}-\mD)^{n+1} .$$
\end{enumerate}
\end{remark}

\begin{example}\label{compute CM}
Let $X$ be a Fano variety of dimension $n$ and $l$ a  natural number. 
We denote $\bP^N:=|-lK_X|$ and $\mD\subset X\times \bP^N$ the universal divisor corresponding to the linear system $|-lK_X|$. For a rational number $0<\epsilon<1$, we compute the CM-line bundle for the family $\pi: (X\times \bP^N, \frac{\epsilon}{l}\mD; \mL)\to \bP^N$, where $\mL\sim_\bQ -K_{X\times \bP^N/\bP^N}-\frac{\epsilon}{l}\mD$.
As the base is of Picard number 1, it suffices to consider a base change for a general rational curve $\bP^1\hookrightarrow \bP^N$, thus we get a family $(X\times \bP^1, \frac{\epsilon}{l}\mD_{\bP^1}; \mL_{\bP^1})\to \bP^1$, still denoted by $\pi$ for convenience. By Remark \ref{remark on CM} we have
$$\lambda_{CM,(X\times \bP^1, \frac{\epsilon}{l}\mD_{\bP^1}, \mL_{\bP^1};\pi )}\sim_\bQ -\pi_*(-K_{X\times \bP^1/\bP^1}-\frac{\epsilon}{l}\mD_{\bP^1})^{n+1}\sim_\bQ \frac{1}{l}(n+1)\epsilon (1-\epsilon)^n(-K_X)^n
\mO_{\bP^1}(1).
$$
Therefore $\lambda_{CM,(X\times \bP^N, \frac{\epsilon}{l}\mD, \mL;\pi )}$ is an ample $\bQ$-line bundle on $\bP^N$. Thus the GIT-stability of $D\in |-lK_X|$ under the action of $\Aut(X)$ with respect to $\lambda_{CM,(X\times \bP^N, \frac{\epsilon}{l}\mD, \mL;\pi )}$ is the same as that with respect to $\mO_{\bP^N}(1)$.
\end{example}

\section{Fano degenerations of a Fano variety}\label{sec:3}

In this section, we fix $X$ to be a Fano variety of dimension $n$.

\begin{definition}
We say that a variety $Y$ is a Fano degeneration of $X$ if there is a $\bQ$- Gorenstein flat family $\mX\to C$ over a smooth pointed curve $0\in C$ such that
\begin{enumerate}
\item $-K_{\mX/C}$ is a relative ample $\bQ$-line bundle,
\item for $t\ne 0$, $\mX_t\cong X$,
\item $\mX_0\cong Y$ is a Fano variety.
\end{enumerate}
\end{definition}

Fix a rational number $0<\epsilon_0<1$, we consider the set $\mF$ of Fano varieties such that $Y\in \mF$ if and only if
\begin{enumerate}
\item $Y$ is a Fano degeneration of $X$,
\item the pair $(Y,cD)$ is K-semistable for some rational $0<c <1-\epsilon_0$ and some $D\sim_\bQ -K_Y$.
\end{enumerate}
Note that the choice of $\epsilon_0\in (0,1)$ does not affect the proof, we just fix such a number.
We have the following lemma.
\begin{lemma}\label{gap lemma}
Notation as above, $\mF$ is contained in a bounded family, and there is a rational number $0<\eta<1$ depending only on $X$ and $\epsilon_0$ such that $Y$ is K-semistable once $\delta(Y)\geq \eta$.
\end{lemma}

\begin{proof}
For $Y\in \mF$, by valuative criterion (see Section \ref{sec value}), we have $\frac{A_{Y,cD}(E)}{S_{Y,cD}(E)}\geq 1 $ for any prime divisor $E$ over $Y$. Note that 
$$S_{Y,cD}(E)=(1-c)S_Y(E),$$ 
thus one sees the following
$$\frac{A_Y(E)}{S_Y(E)}\geq 1-c> \epsilon_0,$$ 
which implies that $\delta(Y)\geq \epsilon_0$. By the relationship between delta invariant and alpha invariant (e.g. \cite{BJ20}), we see that $\alpha(Y)$ also admits a positive lower bound depending only on $\epsilon_0$.  
As $\vol(-K_Y)=\vol(-K_X)$, combining the work \cite{Jiang20}, we know that $\mF$ is contained in a bounded family. By \cite{BLX22,Xu20}, the set 
$$\{\min\{\delta(Y),1\}| Y\in \mF \}$$ 
is finite. Then there exists a rational number $0<\eta<1$ such that $Y$ is K-semistable if $\delta(Y)\geq  \eta$ for any $Y\in \mF$. The proof is finished.
\end{proof}

We are ready to prove the following theorem, which is related to \cite[Theorem 1.2]{ADL19}.
\begin{theorem}\label{capacity}
Notation as above, assume $Y\in \mF$, then there exists a rational number $0<t_1<1$ depending only on $X$ such that if $(Y, tD)$ is K-semistable for some rational $t\in (0, t_1)$ and some $D\sim_\bQ -K_Y$, then $Y$ is K-semistable.
\end{theorem}

\begin{proof}
 Choose a positive rational number $t_1< \min\{1-\epsilon_0,1-\eta\}$, we show that it satisfies our requirement. By Lemma \ref{gap lemma}, it suffices to show that $\delta(Y)\geq \eta$. To see this,  let $E$ be any prime divisor over $Y$, we have the following computation:
$$\frac{A_Y(E)}{S_Y(E)}\geq \frac{A_{Y, tD}(E)}{S_Y(E)}=(1-t)\frac{A_{Y,tD}(E)}{S_{Y,tD}(E)}\geq 1-t\geq \eta. $$ 
Thus we obtain $\delta(Y)\geq \eta$ by valuative criterion. 
\end{proof}

\begin{theorem}\label{key}
Suppose $X$ is a K-polystable Fano variety and $Y$ is a Fano degeneration of $X$. If $Y$ is K-semistable, then $Y\cong X$.
\end{theorem}

\begin{proof}
Let $f: \mX\to C$ be a $\bQ$-Gorenstein flat family over a smooth pointed curve $0_1\in C$ which produces the degeneration from $X$ to $Y$. By \cite[Theorem 1.1]{BX19}, we know there is a degeneration from $Y$ to $X$ via a special test configuration, denoted by $g: \mY\to \bA^1$. We denote the origin $0_2\in \bA^1$, then $\mY_{0_2}\cong X$. Choose a sufficiently large natural number $r$ such that $-rK_\mX$ and $-rK_\mY$ are both very ample, and $\mX$ (resp. $\mY$) can be embedded into $C\times \bP^N$ (resp. $\bA^1\times \bP^N$). Suppose $p(m)=\chi(-mK_X)$, we denote $\Hilb$ the Hilbert scheme whose points parametrize closed sub-varieties in $\bP^N$ with Hilbert polynomial $p(m)$. Then the morphism $f$ (resp. $g$) induces a morphism $C\to \Hilb$ (resp. $\bA^1\to \Hilb$), with $C\setminus 0_1$ and $0_2$ (resp. $\bA^1\setminus 0_2$ and $0_1$) being sent to $[X]\in \Hilb$ (resp. $[Y]\in \Hilb$). Thus we see $\PGL(N+1)[X]\subset \overline{\PGL(N+1)[Y]}$ and $\PGL(N+1)[Y]\subset \overline{\PGL(N+1)[X]}$. Suppose $Y$ is not isomorphic to $X$, then the containments are strict, and we directly have $\dim \PGL(n+1)[Y]<\dim \PGL(n+1)[X]< \dim \PGL(n+1)[Y]$, contradiction.
\end{proof}

\section{Proof of the main result}\label{sec:4}

In this section, we fix $X$ to be a K-polystable Fano variety of dimension $n$. Let $l$ be a positive natural number such that $(X, \frac{1}{l}B)$ is log canonical for some $B\in |-lK_X|$. 

\begin{definition}
We say that the log Fano pair $(Y,\frac{\epsilon}{l}\tD)$ is a degeneration of $(X,\frac{\epsilon}{l}|-lK_X|)$ if there is a $\bQ$-Gorenstein flat family $(\mX,\frac{\epsilon}{l}\mD)\to C$ over a smooth pointed curve $0\in C$ such that 
\begin{enumerate}
\item $-K_\mX-\frac{\epsilon}{l}\mD$ is a relative ample $\bQ$-line bundle,
\item for each $t\ne 0$, $\mX_t\cong X$ and $\mD_t\in |-lK_{\mX_t}|$,
\item $(\mX_0, \mD_0)\cong (Y, \tD)$. 
\end{enumerate}
\end{definition}

Assume $\epsilon<1-\epsilon_0$ (see the definition of $\mF$ in Section \ref{sec:3}). First observe that there exists an Artin stack, denoted by $\mM^K_{X,l,\epsilon}$, generically parametrizing K-semistable log Fano pairs of the form $(X, \frac{\epsilon}{l}D)$, where $D\in |-lK_X|$, and $\mM^K_{X,l,\epsilon}$ admits a good moduli space. To see this point, just note that by \cite[Section 2.6]{XZ20b}, there exists a finite type Artin stack admitting a good moduli space and containing $\mM^K_{X,l,\epsilon}$ as a closed substack with reduced stack structure\footnote{We thank the refree for pointing out this.}.   
It is clear that $\mM^K_{X,l,\epsilon}$  parametrizes all K-semistable log Fano pairs of the form $(Y, \frac{\epsilon}{l}\tD)$, which is the degeneration of $(X,\frac{\epsilon}{l}|-lK_X|)$. Put $\mM^{GIT}_{X,l}:=[(|-lK_X|)^{ss}/\Aut(X)]$. By taking good moduli spaces, we denote $M^K_{X,l,\epsilon}$ (resp. $M^{GIT}_{X,l}$) to be the descent of $\mM^K_{X,l,\epsilon}$ (resp. $\mM^{GIT}_{X,l}$) whose closed points parametrize K-polystable (resp. GIT-polystable) elements.
We first have the following lemma.

\begin{lemma}\label{non-empty}
Notation as above, for general $D\in |-lK_X|$, the pair $(X,\frac{\epsilon}{l}D)$ is K-semisatble for any rational number $0<\epsilon<1$.
\end{lemma}

\begin{proof}
Since $(X, \frac{1}{l}B)$ is log canonical for some $B\in |-lK_X|$, we see that 
for general $D\in |-lK_X|$, the log pair $(X,\frac{1}{l}D)$ is an lc log Calabi-Yau pair, thus being K-semistable by \cite[Corollary 9.4]{BHJ17} or \cite[Theorem 1.5]{Oda13}. By the interpolation property of K-stability (e.g. \cite[Proposition 2.13]{ADL19}), we at once see that $(X,\frac{\epsilon}{l}D)$ is K-semistable for any general $D\in |-lK_X|$ and any rational $0<\epsilon<1$. In particular, the stack $\mM^K_{X,l,\epsilon}$ is not empty.
\end{proof}

We are ready to prove the main result.

\begin{theorem}{\rm (= Theorem \ref{main thm})}\label{thm: main thm}
Let X be a K-polystable Fano variety and $l$ a positive integer such that $(X, \frac{1}{l}B)$ is log canonical for some $B\in |-lK_X|$.  Then there exists a rational number $0<c_1<1$ depending only on  $X$ and $l$ such that the following two statements are equivalent:
\begin{enumerate}
\item $D\in |-lK_X|$ is GIT-(semi/poly)stable under the action of $\Aut(X)$,
\item the pair $(X, \frac{\epsilon}{l}D)$ is K-(semi/poly)stable for any rational $0<\epsilon< c_1$.
\end{enumerate}
As a consequence, we have an isomorphism $\phi_\epsilon: \mM^K_{X,l,\epsilon} \to \mM^{GIT}_{X,l}$ which induces an isomorphism on the moduli spaces $\phi_\epsilon': M^K_{X,l,\epsilon} \to M^{GIT}_{X,l}$ for any rational $0<\epsilon< c_1$.
\end{theorem}

\begin{proof}
We first assume $0<\epsilon\ll 1$. By Lemma \ref{non-empty}, the stack $\mM^K_{X, l, \epsilon}$ is non-empty.
Suppose $(X, \frac{\epsilon}{l}D)$ is K-(semi/poly)stable for $D\in |-lK_X|$, then $D$ is naturally GIT-(semi/poly)stable with respect to an ample CM-line bundle by Remark \ref{remark on CM} and Example \ref{compute CM}.  

Conversely, suppose $D\in |-lK_X|$ is GIT-(semi/poly)stable under $\Aut(X)$-action, by the non-emptiness of $\mM^K_{X, l,\epsilon}$ and openness of K-semistability (e.g. \cite{BLX22, Xu20}),  there is a family $(\mX, \frac{\epsilon}{l}\mD)\to C$ over a smooth pointed curve $0\in C$ such that 
\begin{enumerate}
\item $(\mX_0, \frac{\epsilon}{l}\mD_0)\cong (X, \frac{\epsilon}{l}D)$;
\item $(\mX_t, \frac{\epsilon}{l}\mD_t)$ is K-semistable for any $t\in C\setminus \{0\}$, where $\mX_t\cong X$ and $\mD_t\in |-lK_{\mX_t}|$.
\end{enumerate}
By the properness of K-moduli, up to a finite base change, one could replace the central fiber $(\mX_0, \frac{\epsilon}{l}\mD_0)\cong (X, \frac{\epsilon}{l}D)$ with a K-semistable log Fano pair $(X', \frac{\epsilon}{l}D')$ and we denote the new family by $(\mX', \frac{\epsilon}{l}\mD')\to C$ for convenience. We claim that it suffices to show that $\mX'_0\cong X$. Suppose $\mX_0'\cong X$, then $D'\in |-lK_X|$ and it is GIT-semistable under $\Aut(X)$-action  by what we have proved. By the separatedness of GIT moduli space, we know that $D$ and $D'$ lie on the same orbit under $\Aut(X)$-action. Thus $(X, \frac{\epsilon}{l}D)$ is K-(semi/poly)stable.

Next we show $\mX_0'\cong X$. Since $\epsilon$ is sufficiently small, by Theorem \ref{capacity}, we know that $\mX_0'$ is K-semistable. Note that $X$ is K-polystable, then $\mX_0'\cong X$ by Theorem \ref{key}.

We are ready to establish the isomorphism between stacks. For any $\bC$-valued point $[(Y, \frac{\epsilon}{l}\tilde{D})]\in \mM^K_{X, l, \epsilon}$, we see that $Y\cong X$ and $\tilde{D}\in \frac{1}{l}|-lK_X|$ is GIT-semistable.
Thus  we have the following morphism
$$\phi_\epsilon: \mM^K_{X, l,\epsilon}\to \mM^{GIT}_{X,l},\ \ \ [(X, \frac{\epsilon}{l}D)]\mapsto [D].$$
Similar to the argument at the end of the proof of \cite[Theorem 5.2]{ADL19}, $\phi_\epsilon$ is an isomorphism. Thus it descends to an isomorphism between good moduli spaces in the sense of Alper (e.g. \cite{Alper13}), i.e., $M^K_{X, l,\epsilon}\to M^{GIT}_{X,l}$
is an isomorphism.

The existence of $c_1$ is given by \cite{Zhou21},   where we show that $\mM^K_{X, l, \epsilon}$ (resp. $M_{X, l, \epsilon}^K$) does not change as we vary $\epsilon$ in $(0, c_1)$, and $c_1$ depends on $X$ and the coefficient data $l$.
The proof is finished.
\end{proof}

\bibliography{reference.bib}
\end{document}